\documentclass[11pt,a4paper, twoside]{article}
\usepackage{amsmath,amssymb,amsthm,amsfonts,mathrsfs,amscd,environ,bbm}
\usepackage{latexsym,enumerate,color,geometry}
\usepackage{xcolor}
\usepackage{verbatim}
 \NewEnviron{ews}{%
\begin{equation}\begin{split}
  \BODY
\end{split}\end{equation}
}

\NewEnviron{ews*}{%
\begin{equation*}\begin{split}
  \BODY
\end{split}\end{equation*}
}

\def\beg{\begin}
\def\bequ{\begin{equation}}
\def\enqu{\end{equation}}
\def\bes{\begin{split}}
\def\ens{\end{split}}
\def\bews{\begin{ews}}
\def\beqn{\begin{eqnarray}}
\def\enqn{\end{eqnarray}}
\def\beq*{\begin{equation*}}
\def\enq*{\end{equation*}}
\def\bqn*{\begin{eqnarray*}}
\def\eqn*{\end{eqnarray*}}
\def\bary{\begin{array}}
\def\eary{\end{array}}
\def\bpma{\begin{pmatrix}}
\def\epma{\end{pmatrix}}
\def\bvma{\begin{Vmatrix}}
\def\evma{\end{Vmatrix}}

 \numberwithin{equation}{section}

\def\al{\alpha}
\def\be{\beta}
\def\ga{\gamma}
\def\de{\delta}

\def\ze{\zeta}
\def\et{\eta}

\def\la{\lambda}

\def\si{\sigma}
\def\ta{\tau}

\def\om{\omega}

\def\Ga{\Gamma}

\def\De{\Delta}

\def\Ph{\Phi}

\def\Om{\Omega}

\def\R{\mathbb R}
\def\P{\mathbb P}
\def\E{\mathbb E}

\def\N{\mathbb N}

\def\SS{\mathbb S}

\def\sF{\mathscr F}

\def\sB{\mathscr B}

\def\sG{\mathscr G}

\def\d{\mathrm{d}}

\def\ff{\frac}
\def\ra{\rightarrow}

\def\<{\langle}
\def\>{\rangle}

\def\tld{\tilde}
\def\we{\wedge}
\def\1{\mathbbm{1}}

\allowdisplaybreaks

\title{{\bf Regime-switching diffusion processes: strong solutions and strong Feller property}
}

\author{
{\bf Shao-Qin Zhang }\\
\footnotesize{School of Statistics and Mathematics, Central University of Finance and Economics, Beijing 100081, China}\\
\footnotesize{Email: zhangsq@mail.bnu.edu.cn}\\
}

\begin{document}
\maketitle

\begin{abstract}
We investigate the existence and uniqueness of strong solutions up to an explosion time for regime-switching diffusion processes in an infinite state space.  Instead of concrete conditions on coefficients, our existence and uniqueness result is established under the general assumption that the diffusion in every fixed environment has a unique non-explosive strong solution. Moreover, non-explosion conditions for regime-switching diffusion processes are given. The strong Feller property is proved by further assuming that the diffusion in every fixed environment generates a strong Feller semigroup. 
\end{abstract}\noindent

\noindent


\vskip 2cm

\section{Introduction}
Many uncertain hybrid systems from financial engineering, wireless communications, biology and etc. can be modeled by the regime-switching diffusion processes(RSDPs for short), see \cite{YinZB} and references.  In these models, the continuous dynamics and discrete events coexist. Various properties of these processes have been concerned recently. For example, \cite{PinP,PinSch,Shao14a,Shao14b,ShaoX13,ShaoX14,YinZ09} had studied the recurrent properties, \cite{GAM,MaoY,ShaoX14,XiZ,YinZB} concerned stability and optimal control, and \cite{Skoro}  investigated the asymptotic properties.  We mention monographs \cite{MaoY} and \cite{YinZB} for system surveys on the research of regime-switching processes.  

Let $\SS$ be the set of all the natural numbers, i.e. $\SS=\{1,2,3,\cdots,n,\cdots\}$, $(\Om,\sF,\P)$ be a completed probability space. A regime-switching diffusion process is a two components process $\{(X_t,\Lambda_t)\}_{t\geq 0}$ described by 
\bequ\label{equ_b1}
\d X_t=b(X_t,\Lambda_t,t)\d t+\si(X_t,\Lambda_t,t)\d W_t,
\enqu
and
\bequ\label{equ_b2}
\P\left(\Lambda_{t+\De t}=j~|~\Lambda_t=i,~(X_s,\Lambda_s),s\leq t\right)
=\left\{\beg{array}{lr}
q_{ij}(X_t)\De t+o(\De t),& i\neq j,\\
1+q_{ii}(X_t)\De t+o(\De t),& i=j,
\end{array}
\right.
\enqu
where $\{W_t\}_{t\geq 0}$ is a Brownian motion on $\R^d$ w.r.t a completed reference $\{\sF_t\}_{t\geq 0}$ and 
$$b:\R^d\times \SS\times [0,\infty)\ra \R^d,~~\si:\R^d\times \SS\times [0,\infty)\ra \R^d\otimes\R^d,~~q_{ij}:\R^d\ra \R,$$
are  measurable functions and 
$$q_{ij}(x)\geq 0,~i\neq j,~~~~~\sum_{j\neq i}q_{ij}(x)\leq -q_{ii}(x).$$
The matrix $Q(x)=(q_{ij}(x))$ is called $Q$-matrix, see \cite{ChenB}. If $Q(x)$ is independent of $x$ and $\Lambda_t$ is independent of $\{W_t\}_{t\geq 0}$, $\{(X_t,\Lambda_t)\}_{t\geq 0}$ is called a state-independent RSDP, otherwise called a state-dependence RSDP. Let 
$$q_i(x)=\sum_{j\neq i}q_{ij}(x),~i\in\SS.$$ If $Q(x)$ is conservative, i.e. $-q_{ii}(x)= q_i(x)$, $x\in\R^d$, then as in \cite[Chapter II-2.1]{Skoro} or \cite{GAM,YinZ09,Shao15}, we can represent $\{(X_t,\Lambda_t)\}_{t\geq 0}$ in the form of a system of stochastic differential equations(SDEs for short) driven by $\{W_t\}$ and a Poisson random measure. Precisely, for each $x\in\R^d$, $\{\Ga_{ij}(x):~i,~j\in\SS\}$ is a family of disjoint  intervals on $[0,\infty)$ constructed as follows
\beg{align*}
&\Ga_{12}(x)=[0,q_{12}(x)),~\Ga_{13}(x)=[q_{12}(x),q_{12}(x)+q_{13}(x)),\cdots\\
&\Ga_{21}(x)=[q_1(x),q_1(x)+q_{21}(x)),~\Ga_{23}(x)=[q_1(x)+q_{21}(x),q_1(x)+q_{21}(x)+q_{23}(x)),\cdots\\
&\Ga_{31}(x)=[q_1(x)+q_2(x),q_1(x)+q_{2}(x)+q_{13}(x)),\cdots\\
&\qquad\cdots
\end{align*}
We set $\Ga_{ii}(x)=\emptyset$ and $\Ga_{ij}(x)=\emptyset$ if $q_{ij}(x)=0$ for $i\neq j$. Define a function $h:\R^d\times\SS\times [0,\infty)\ra\R$  as follows
\beq*
h(x,i,z)=\sum_{j\in\SS}(j-i)\1_{\Ga_{ij}(x)}(z)=\left\{\beg{array}{cl}
j-i,& \mbox{if~} z \in \Ga_{ij}(x),\\
0,& \mbox{otherwise}.
\end{array}
\right.
\enq*
Let $N(\d z,\d t)$ be a Poisson random measure with intensity $\d t\d z$ and independent of the Brownian motion $\{W_t\}_{t\geq 0}$. Then we turn to concern the following equation
\bequ\label{main_equ}
\beg{cases}
\d X_t = b(X_t,\Lambda_t,t)\d t+\si(X_t,\Lambda_t,t)\d W_t,\\
\d \Lambda_t = \int_{0}^\infty h(X_{t-},\Lambda_{t-},z)N(\d z,\d t).
\end{cases}
\enqu
Denote the diffusion process in the fixed environment $i\in\SS$ by $x_t^i$, see \eqref{equ_i} for more precise definition. Showing as in \eqref{main_equ} or \eqref{equ_b1} and \eqref{equ_b2}, one can image that the diffusion component $\{X_t\}_{t\geq 0}$ is a hybrid process of  $\{x^i_t\}_{i\in\SS}$ via jump process $\Lambda_t$. A quite nature question is that if all the diffusion processes $\{\{x_t^i\}_{t\geq 0}\}_{i\in\SS}$ process some comment property(such as pathwise uniqueness, strong Feller property, recurrence, ergodicity and etc.), then what about $\{(X_t,\Lambda_t)\}_{t\geq 0}$? In fact, many works, such as \cite{Shao14b,Shao15,ShaoX14,YinZ09}, give us some insight on this problem. 

There are some papers focusing on existence and uniqueness of strong solutions and strong Feller property for RSDPs. For instance, we refer to \cite{XiY11} for existence and uniqueness of solution in a finite state space, to \cite{Shao15,XiY13,YinZ09} for strong Feller property. It should be pointed out that, \cite{Shao15} established the existence and uniqueness and strong Feller property to the strong solution of state-dependent RSDPs in an infinite state space with some type of non-Lipschitz coefficients. However, existence and uniqueness theorem in \cite{Shao15} can not be applied to more irregular case, such as H\"older continuous coefficients or Sobolev coefficients, which have been intensively studied recently. Here, we only mention \cite{XZ} for a brief survey in their introduction on strong uniqueness of SDEs with singular coefficients  and references therein. Hence one aim of the this paper is to investigate the existence and uniqueness of the regime-switching processes that the coefficients may be singular. Another aim is to concern strong Feller property  of $\{(X_t,\Lambda_t)\}_{t\geq 0}$ under the assumption that the diffusion in every fixed environment processes the same property.  The method used here is different from \cite{Shao15}. Our argument is based on the following simple observation(formally), that is if $\{q_i\}_{i\in\SS}$ are locally bounded and $\sup_{s\in[0,t]}\left(|X_{s-}|+\Lambda_{s-}\right)\leq M$($M\in\SS$), then
$$\bigcup_{s\in[0,t]}{\rm{supp}}\left(h(X_{s-},\Lambda_{s-},\cdot)\right)\subset \bigcup_{|x|\leq M, i\leq M,j\geq 1}\Ga_{i,j}(x)\subset \left[0,\sup_{|x|\leq M} \sum_{k=1}^{M+1}q_k(x)\right],$$
and  $\sup_{|x|\leq M} \sum_{k=1}^{M+1}q_k(x)<\infty$. Let $K=[0,\sup_{|x|\leq M} \sum_{k=1}^{M+1}q_k(x)]$. Then
$$\Lambda_t=\Lambda_0+\int_0^t\int_{0}^\infty h(X_{s-},\Lambda_{s-},z)N(\d z,\d s)=\Lambda_0+\int_0^t\int_{K} h(X_{s-},\Lambda_{s-},z)N(\d z,\d s).$$
But on $K$, $N(K,t):=\int_0^t\int_KN(\d z,\d s)$ is a Poisson process which jumps finite times on any finite  interval. This allows us to generalize the method used in \cite{XiY11}(see also \cite[Theorem IV.9.1]{IW89}), and give an interlacing construction to the solution of \eqref{main_equ} before $(X_t,\Lambda_t)$ escapes a bounded domain. The strong Feller property is also studied following this idea which is also different from \cite{Shao15,YinZ09}, where they got the strong Feller property by a finite-state approximation argument. Since we do not assume the concrete conditions on coefficients of the  diffusion part,  our results can be applied to the case that the diffusion component may be irregular and degenerated.  

The rest part of the paper is organized as follows. In Section 2, we discuss the existence, uniqueness and non-explosion of the solution to \eqref{main_equ}. In Section 3, we study the strong Feller property.

\section{Existence, uniqueness and non-explosion}
Throughout this paper, we assume that the $Q$-matrix $Q(x)$ is conservative, that is
$$-q_{ii}(x)= q_i(x)\left(=\sum_{j\neq i}q_{ij}(x)\right),$$
for more details on continuous time Markov chains, one can consult \cite{ChenB} .  By this assumption, then we shall prove the existence of a unique strong solution to the equation \eqref{main_equ}.

We introduce the following hypothesis
\begin{description}
  \item [{\bf Hypothesis}\  {\bf (H)}] For each $i\in\SS$, $t_0\geq 0$, the following equation
\bequ\label{equ_i}
\d x_t=b(x_t,i,t+t_0)\d t+\si(x_t,i,t+t_0)\d W_t
\enqu
has a unique non-explosive strong solution.
\end{description}
Let
$$W^d=C([0,\infty),\R^d),~~~W_{0}^d=\left\{\om\in W^d~|~\om(0)=0\right\},$$
with the topology of locally uniformly convergence on $[0,\infty)$, $\sB(W^d)$ be the Borel $\si$-field of $W^d$ and $\sB_t(W^d)$ be the sub-$\si$-field of $\sB(W^d)$ generated by $\om(s),~0\leq s\leq t$. We recall that, for each $(i,t_0)\in\SS\times[0,\infty)$, a $d$ dimensional continuous process $x^{i,t_0}:=\{x^{i,t_0}_t\}_{t\geq 0}$ defined on a probability space with a reference family and a Brownian motion $\{W_t\}_{t\geq 0}$ is a strong solution of \eqref{equ_i}, if there is a function $F(i,t_0,\cdot,\cdot):\R^d\times W_0^d\ra W^d$ such that
\beg{enumerate}
\item [$(1)$] for all  Borel probability measure $\mu$ on $\R^d$, there is a function $\tld F_\mu(i,t_0,\cdot,\cdot):\R^d\times W_0^d\ra W^d$  which is $\overline{\sB(\R^d\times W^d_0)}^{\mu\times\P^W}/\sB(W^d)$-measurable, and for $\mu$-a.e. $x\in\R^d$ it holds that
$$F(i,t_0,x,\om)=\tld F_\mu(i,t_0,x,\om),~\P^W\mbox{-a.s.~} \om\in W^d_0,$$
where $\P^W$ is the Wiener measure on $W_0^d$,
\item [$(2)$] for each $x\in\R^d$, $F(i,t_0,x,\cdot):W_0^d\ra W^d$ is $\overline{\sB_t(W_0^d)}^{\P^W}/\sB_t(W^d)$ measurable for each $t\geq 0$, and $x^{i,t_0}=F(i,t_0,x^{i,t_0}_0,W_\cdot)$.
\end{enumerate}
For more details on the strong solution of SDEs, we refer to \cite[Chapter IV. Definition 1.6]{IW89}. In the following, we shall denote the solution of \eqref{equ_i} by $x_t^{x,i,t_0}$ if it starts from $x\in\R^d$ with parameters $i,t_0$ in the coefficients, and omit $t_0$ in superscript if $t_0=0$.  Here, we emphasize that in $(2)$, the initial value of $x^{i,t_0}_t$, i.e. $x_0^{i,t_0}$, can be chosen arbitrary though we write parameters $i,t_0$ in superscript.

In general, the solution to \eqref{main_equ} can be explosive. So we consider the local solution.

\beg{defn}
An adapted process $\{(X_t,\Lambda_t)\}_{t\in[0,\ta)}$ with the first component is continuous and the second one is cadlag is called a local solution of \eqref{main_equ} with life time $\ta$, if $\ta>0$ is a stopping time such that $\P$-a.s. $\limsup_{t\uparrow \ta}\left(|X_t|+\Lambda_t\right)=\infty$ holds on $\{\ta<\infty\}$ and $\P$-a.s.
 \beq*
\beg{cases}
X_t =X_0+ \int_0^t b(X_s,\Lambda_s,s)\d s+\int_0^t\si(X_s,\Lambda_s,s)\d W_s,\\
\Lambda_t =\Lambda_0 +\int_0^t\int_{0}^\infty h(X_{s},\Lambda_{s-},z)N(\d z,\d s),~t\in[0,\ta).
\end{cases}
\enq*
\end{defn}
Our first result is the following
\beg{thm}\label{path_unique}
Assume that {\bf(H)} holds, and for  each $k\in\SS$, $q_k$ is locally bounded. Then for each initial value $(x,i)\in\R^d\times\SS$, the equation \eqref{main_equ} has a unique local strong solution $\{(X_t,\Lambda_t)\}_{t\in[0,\ta)}$ with life time $\ta$.
\end{thm}

We can extend Theorem \ref{path_unique} to some cases that the diffusion in each environment may be explosive. 
\beg{cor}\label{cor_loc}
For all $M>0$, let
$$b_M(x,i,t)=b(x,i,t)\1_{\{|x|\leq M,t\leq M\}}(x),~\si_M(x,i,t)=\si(x,i,t)\1_{\{|x|\leq M,t\leq M\}}(x).$$
If for all $M>0$,~$i\in\SS$, the equation
$$\d x_t^M=b_M(x_t^M,i,t)\d t+\si_M(x_t^M,i,t)\d W_t$$
has a unique non-explosive strong solution, and $q_i(\cdot)$ is locally bounded, then the equation \eqref{main_equ} has a unique local strong solution $\{(X_t,\Lambda_t)\}_{t\in[0,\ta)}$ with lift time $\ta$.
\end{cor}

Before the proof of our results, we present a lemma on a SDE(see \eqref{equ_random})  with  random coefficients related to \eqref{equ_i}. Let  $\left(\tld \P,\tld\Om,\tld\sF\right)$ be some completed probability space, $\left\{\tld W_t\right\}_{t\geq 0}$ be a Brownian motion on this probability space w.r.t. some completed reference $\left\{\tld\sF_t\right\}_{t\geq 0}$. 
Let $(\et,\ze,\xi):\tld\Om\ra  \SS\times (0,\infty)\times \R^d$ be a random vector which is independent of  $\{\tld W_t\}_{t\geq 0}$ and $(\et,\ze,\xi)\in\tld \sF_0$. Denote the distribution of $(\et,\ze,\xi)$ by $\mu$. Let 
$$\mathring{\sG}_t=\si\left(\left\{\et,\ze,\xi,\tld W_s~|~s\leq t\right\}\right),$$
$\sG_t$ be the completion of $\mathring{\sG}_t$ by $\tld\P$.
Consider the following SDEs with random coefficients
\bequ\label{equ_random}
\d x_t=b(x_t,\et,t+\ze)\d t+\si(x_t,\et,t+\ze)\d \tld W_t,~x_0=\xi.
\enqu
We say a continuous adapted process $\{x_t\}_{t\geq 0}$ is a solution to \eqref{equ_random} with a Brownian motion $\{\tld W_t\}_{t\geq 0}$, if 
$$\int^{t}_0\left|b(x_s,\et,s+\ze)\right|+\left|\si(x_s,\et,s+\ze)\right|^2\d s<\infty,~t\geq 0,~\tld\P\mbox{-a.s.}$$
and
$$x_{t}=x_0+\int_0^{t} b(x_s,\et,s+\ze)\d s+\int_0^{t}\si(x_s,\et,s+\ze)\d \tld W_s,~t\geq 0,~\tld\P\mbox{-a.s.}$$
Due to the hypothesis {\bf{(H)}}, there exists $F:\SS\times (0,\infty)\times\R^d\times W^d_{0}\ra W^d$ such that $\{F(i,t_0,x,\tld W_\cdot)_t\}_{t\geq 0}$ is the solution of \eqref{equ_i} starting from $x$. Since $(\et,\ze,\xi)$ is independent of $\{\tld W\}_{t\geq 0}$ and $(\xi,\et,\ze)\in\tld \sF_0$, $\{F(\xi,\et,\ze,\tld W_\cdot)_t\}_{t\geq 0}$ is a solution of \eqref{equ_random}.
\beg{lem}\label{lem_path}
Assume that {\bf(H)} holds. Let $\{x'_t\}_{t\geq 0}$ be a solution of \eqref{equ_random}. Then $x'_t=F(\et,\ze,\xi,\tld W_\cdot)_t$, $t\geq 0$,~$\tld\P$-a.s.  
\end{lem}
\beg{proof}
It is clear that $x'_\cdot\in\sG_{\infty}$ and $x'_t\in\sG_t$, so there is $F'_\mu:~\SS\times (0,\infty)\times\R^d\times  W^d_{0}\ra W^d$ such that 
$$F'_\mu(\et,\ze,\xi,\tld W_{\cdot})_t\in\tld \sF_t,~~~\mbox{and}~~~\tld\P\left(x'=F'_\mu(\et,\ze,\xi,\tld W_\cdot)\right)=1.$$ 
Due to \cite[Lemma V.10.1]{RogersW}, there is $\Ph:\SS\times (0,\infty)\times \R^d\times W^d_{0}\ra W^d$ such that
$$\int_0^t \si(F'_\mu(\et,\ze,\xi,\tld W_{\cdot})_s,\et,s+\ze)\d \tld W_s=\Ph(\et,\ze,\xi,\tld W_{\cdot})_t,~t\geq 0.$$
Then we have $\tld \P\mbox{-a.s.}$
$$F'_\mu(\et,\ze,\xi,\tld W_{\cdot})_t=\xi+\Ph(\et,\ze,\xi,\tld W_{\cdot})_t+\int_0^tb(F'_\mu(\et,\ze,\xi,\tld W_{\cdot})_s,\et,s+\ze)\d s,~t\geq 0.$$
Hence, for $\mu$-a.s. $(i,t_0,x)\in\SS\times (0,\infty)\times\R^d$,  $\tld \P\mbox{-a.s.}$
$$F'_\mu(i,t_0,x,\tld W_{\cdot})_t=x+\Ph(i,t_0,x,\tld W_{\cdot})_t+\int_0^tb(F'_\mu(i,t_0,x,\tld W_{\cdot})_s,i,s+t_0)\d s,~t\geq 0.$$
That means $F'_\mu(i,t_0,x,\tld W_{\cdot})_t$ is a solution of \eqref{equ_i} with $W$ replaced by $\tld W$. If 
$$\tld \P\left(x'\neq F(\et,\ze,\xi,\tld W_\cdot)\right)>0,$$
then
$$\tld\P\left(F(\et,\ze,\xi,\tld W_\cdot)\neq F'_\mu(\et,\ze,\xi,\tld W_\cdot)\right)>0.$$
Since $(\et,\ze,\xi)$ is independent of $\{\tld W_t\}_{t\geq 0}$, there is $A\subset\SS\times(0,\infty)\times\R^d$ with $\mu(A)>0$ such that
$$\tld\P\left(F(i,t_0,x,\tld W_\cdot)\neq F'_\mu(i,t_0,x,\tld W_\cdot)\right)>0,~~(i,t_0,x)\in A.$$
It contradicts  {(\bf{H})},  since $\{F(i,t_0,x,\tld W_\cdot)_t\}_{t\geq 0}$ and $\{F'_\mu(i,t_0,x, \tld W_\cdot)_t\}_{t\geq 0}$ are two solutions of \eqref{equ_i} on the probability space $\left(\tld\Om, \tld\sF,\tld \P\right)$ with the same Brownian motion $\tld W$.

\end{proof}

\medskip
{\noindent\textbf{\emph{Proof of Theorem \ref{path_unique}}}~}\\
Let $K>0$. We also denote $K=[0,K]$ without confusion. We first show the existence and uniqueness of  solutions to the following equation
\bequ\label{equ_K}
\beg{cases}
\d X^K_t = b(X^K_t,\Lambda_t,t)\d t+\si(X^K_t,\Lambda^K_t,t)\d W_t,~~X^K_0=x\\
\d \Lambda^K_t = \int_{K} h(X^K_{t-},\Lambda^K_{t-},z)N(\d z,\d t),~~\Lambda^K_0=i.
\end{cases}
\enqu
Let 
$$N(K,t)=\int_0^t\int_{K}N(\d z,\d s),~~N_t^K=\int_0^t\int_KzN(\d z,\d t).$$
Then $\{N(K,t)\}_{t\geq 0}$ is a Poisson process and $\{N_t^K\}_{t\geq 0}$ is a compound Poisson process. 
Let
$$T_1^K=\inf\{t>0~|~N(K,t)>0\},~~T_{n+1}^K=\inf\{t>T_{n}^K~|~N(K,t)>n\},~n\geq 1.$$
Denote the size of the $n$-th jump of $N^K_t$ by $\De N_n^K$. Then $\{T_n^K\}_{n\in\N}$, $\{\De N_n^K\}_{n\in\N}$ and $\{W_t\}_{t\geq 0}$ are independent of each other and
$$0<T_1^K<T_2^K<\cdots<T_n^K<\cdots<\infty.$$

Set $(X_t^K,\Lambda_t^K)=(x_t^{x,i,0},i)$, $t\in[0,T_1^K)$, and $X_{T_1^K}^K=x_{T_1^K}^{x,i,0}$. 
Hence $X_t^K$ satisfies that
\beg{align*}
X_{t\we T_1^K}^K&=x+\int_0^{t\we T_1^K}b(x_t^{x,i,0},i,s)\d t+\int_0^{t\we T_1^K}\si(x_t^{x,i,0},i,s)\d W_s\\
&=x+\int_0^{t\we T_1^K}b(X_s^K,\Lambda_s^K,s)\d s+\int_0^{t\we T_1^K}\si(X_s^K,\Lambda_s^K,s)\d W_s.
\end{align*}
Then at $T_1^K$, we have
$$\Lambda_{T_1^K}=i+h(X^K_{T_1^K-},i,\De N_1^K)=i+h(x^{x,i,0}_{T_1^K-},i,\De N_1^K).$$
So, according to {\bf(H)},  the solution to \eqref{equ_K} on $[0,T_1^K]$ exists uniquely and $(X_t^K,\Lambda_t^K)$ can be constructed as above.

Next, set $X^{1,K}_0=X^K_{T_1}$, $\Lambda^{1,K}_0=\Lambda_{T_1^K}$ and
$$W^1_t=W_{t+T_1^K}-W_{T_1^K},~~N^1_t=N^K_{T_1^K+t}-N^K_{T_1^K},~~t\geq 0.$$
According to the strong Markov property, $\{W^1_t\}_{t\geq 0}$ is also a Brownian motion and $N^1_t$ is compound Poisson process w.r.t $\{\sF_{T_1^K+t}\}_{t\geq 0}$. Let $N^1(\d z,\d t)$ be the  corresponding Poisson random measure of $N^1_t$.
We consider 
\bequ\label{equ_shift}
\beg{cases}
X^{1,K}_t =X^K_{T_1}+\int_0^t b(X^{1,K}_s,\Lambda_s,s+T_1^K)\d s+\si(X^{1,K}_s,\Lambda^{1,K}_s,s+T_1^K)\d W^1_s,\\
\Lambda^{1,K}_t =\Lambda_{T_1^K}+ \int_0^t\int_{K} h(X^{1,K}_{s-},\Lambda^{1,K}_{s-},z)N^1(\d z,\d s).
\end{cases}
\enqu
Denote the time of the first arrival of $N^1(K,t)$ by $T^{1,K}_1$, then $T^{1,K}_1$ coincide with $T_2^K-T_1^K$.
It is clear that  we can set
$$\Lambda_t^{1,K}=\Lambda^{1,K}_0,~~t<T_1^{1,K}.$$
Then $X_{t}^{1,K}$ should satisfy 
$$X_{t}^{1,K}=X_0^{1,K}+\int_0^{t}b(X_s^{1,K},\Lambda_0^{1,K},s+T_1^K)\d s+\int_0^{t}\si(X_s^{1,K},\Lambda_0^{1,K},s+T_1^K)\d W^1_s,~t<T_1^{1,K},$$
with $\{W^1_s\}_{s\geq 0}$ is independent of $\sF_{T_1^K}$ and $\Lambda_0^{1,K},T_1^K\in\sF_{T_1^K}$. So, according to  Lemma \ref{lem_path},  \eqref{equ_shift} has a unique solution on $[0,T_1^{1,K}]$ and we can construct it as follows
$$X_t^{1,K}:=F(X_0^{1,K},\Lambda_0^{1,K}, T_1^K,W^1_\cdot)_t,~~0<t<T_1^{1,K}.$$
Then 
$$\Lambda_{T_1^{1,K}}=\Lambda_0^{1,K}+h(X^{1,K}_{T_1^{1,K}-},\Lambda_0^{1,K},\De N_1^K).$$
Set $X_{T_1^{1,K}}^{1,K}=F(X_0^{1,K},\Lambda_0^{1,K}, T_1^K,W^1_\cdot)_{T_1^{1,K}}$. Define
$$\left(X_{t}^K,\Lambda^K_t\right)=\left(X_{t-T_1^K}^{1,K},\Lambda_{t-T_1^K}^{1,K}\right),~~t\in[T_1^K,~T_2^K].$$
By continuing this process, we can determine the solution of \eqref{equ_K} on all $[0,T_n^K]$, $n\geq 1$. Since $T_n^K\ra \infty$ as  $n\ra\infty$ and Hypothesis {\bf(H)}, we obtain a unique non-explosive solution  $(X_{t}^K,\Lambda^K_t)$ to \eqref{equ_K}.

Next, let $M,K\in \SS$ such that
\bequ\label{KM}
M>|x|+i,~K\geq \sup_{|y|\leq M}\sum_{k=1}^{M+1}q_k(y).
\enqu
Set 
$$\ta_M^K=\inf\{t\geq 0~|~|X_t^K|+\Lambda_t^K\geq M\},~~\ta_M^{K+1}=\inf\{t\geq 0~|~|X_t^{K+1}|+\Lambda_t^{K+1}\geq M\},$$
where $\left(X_t^{K+1},\Lambda_t^{K+1}\right)$ is the solution of \eqref{equ_K} with $K$ replaced by $K+1$. Since $T_1^K>0$, $\P$-a.s. and 
$$\left(X_t^K,\Lambda_t^K\right)=(x_t^{x_0,i},i),~t<T_1^K,$$
we have  $\ta_M^K>0$, $\P$-a.s. Similarly, $\ta_M^{K+1}>0$, $\P$-a.s. Note that
\beg{align*}
{\rm{supp}}\left(h(X^{K+1}_{t-},\Lambda^{K+1}_{t-},\cdot)\right)&\subset \bigcup_{|x|\leq M, i\leq M,j\geq 1}\Ga_{i,j}(x)\\
&\subset \left[0,\sup_{|x|\leq M} \sum_{k=1}^{M+1}q_k(x)\right]\subset [0,K],~~t\leq\ta^{K+1}_M.
\end{align*}
Then
\beg{align*}
X^{K+1}_t&=x+\int_0^{t}b(X^{K+1}_s,\Lambda^{K+1}_s,s)\d s+\int_0^{t}\si(X^{K+1}_s,\Lambda^{K+1}_s,s)\d W_s,\\
\Lambda^{K+1}_t&=i+\int_0^t\int_{K+1} h(X^{K+1}_{s-},\Lambda^{K+1}_{s-},z)N(\d z,\d s)\nonumber\\
&=i+\int_0^t\int_{K}h(X^{K+1}_{s-},\Lambda^{K+1}_{s-},z)N(\d z,\d s),~~t\leq\ta^{K+1}_M.
\end{align*}
Due to the uniqueness of \eqref{equ_K}, we get that  $\ta_M^{K+1}\leq \ta_M^K$ and
$$\left(X_t^{K+1},\Lambda_t^{K+1}\right)=\left(X_t^K,\Lambda_t^K\right),~t\leq\ta_M^{K+1}.$$
However, 
$$|X_{\ta_M^{K+1}}^{K+1}|+\Lambda^{K+1}_{\ta_M^{K+1}}\geq M,$$
so, in fact, we have $\ta_M^{K+1}=\ta_M^{K}$. Hence, we set 
$$\ta_M=\lim_{K\ra\infty}\ta_M^K,~~\ta=\lim_{M\ra\infty}\ta_M.$$
$\ta_M$ is an increasing sequence of stopping time and $\ta_M>0$. Now, for each $M\in\SS$ with $M>|x|+i$ and $K$ large enough(recalling that \eqref{KM}), we define
$$\left(X_t,\Lambda_t\right)=\left(X_t^K,\Lambda_t^K\right),~t<\ta_M=\ta_M^K.$$
Then $(X_t,\Lambda_t)$ satisfies \eqref{main_equ} for $t<\ta$ and 
$$\limsup_{t\uparrow\ta}\left(|X_t|+\Lambda_t\right)=\infty.$$

Lastly, we prove the pathwise uniqueness of solutions to \eqref{main_equ}. Let $(X_t,\Lambda_t)$ and $(\tld X_t,\tld \Lambda_t)$ be two solutions in the same probability space with the same Brownian motion and Poisson random measure. For each $M\in\SS$ with $M>|x|+i$, let
\bequ\label{ta_M}
\ta_M=\inf\left\{t\geq 0~\Big|~|X_t|+|\Lambda_t|\geq M\right\},~~\tld\ta_M=\inf\left\{t\geq 0~\Big|~|\tld X_t|+|\tld\Lambda_t|\geq M\right\}.
\enqu
Let $K= \left[0,\sup_{|x|\leq M} \sum_{k=1}^{M+1}q_k(x)\right]$. Then
\bequ\label{equ_MK}
\beg{cases}
X_t=x+\int_0^{t}b(X_s,\Lambda_s,s)\d s+\int_0^{t}\si(X_s,\Lambda_s,s)\d W_s,\\
\Lambda_t=i+\int_0^t\int_{K}h(X_{s-},\Lambda_{s-},z)N(\d z,\d s),~~t\leq\ta_M.
\end{cases}
\enqu
That means  $\left(X_t,\Lambda_t\right)$ is a solution of \eqref{equ_K} before $\ta_M$. So does for $\left(\tld X_t,\tld \Lambda_t\right)$. Then by uniqueness of \eqref{equ_K}, we obtain that
$$\left(X_{t\we\ta_M},\Lambda_{t\we\ta_M}\right)=\left(\tld X_{t\we\tld\ta_M},\tld \Lambda_{t\we\tld\ta_M}\right).$$
Therefore we get the pathwise uniqueness of solutions to \eqref{main_equ}.

\qed

\medskip
{\noindent\textbf{\emph{Proof of Corollary \ref{cor_loc}}}~}\\
We adapt the stopping time skill used in \cite[Theorem 1.1]{ZhangX}. Let $M\in\SS$ such that $M>|x|+i$. We first consider the following equation
\bequ\label{equ_main_M}
\beg{cases}
\d X^M_t = b_M(X^M_t,\Lambda^M_t,t)\d t+\si_M(X^M_t,\Lambda^M_t,t)\d W_t,X_0^M=x,\\
\d \Lambda^M_t = \int_{0}^\infty h(X^M_{t-},\Lambda^M_{t-},z)N(\d z,\d t),\Lambda_0^M=i.
\end{cases}
\enqu
By Theorem \ref{path_unique} and our assumption, \eqref{equ_main_M} has a unique local strong solution. Let $(X^{M+1}_t,\Lambda_t^{M+1})$ be the solution of \eqref{equ_main_M} with $M$ replaced by $M+1$,
\beg{align*}
\ta_M^M &=\inf\{t\in[0,M)~|~|X_t^M|+\Lambda_t^M\geq M\},\\
\ta_{M}^{M+1}&=\inf\{t\in[0,M)~|~|X_t^{M+1}|+\Lambda_t^{M+1}\geq M\}.
\end{align*}
Since $|X_{t\we\ta_M^{M+1}}|\leq M$, we have that
\beg{align*}
X_{t\we\ta_M^{M+1}}^{M+1}&=x+\int_0^{t\we\ta_{M}^{M+1}}b_{M+1}(X_s^{M+1},\Lambda_s^{M+1},s)\d s\\
&\qquad+\int_0^{t\we\ta_{M}^{M+1}}\si_{M+1}(X_s^{M+1},\Lambda_s^{M+1},s)\d W_s\\
&=x+\int_0^{t\we\ta_{M}^{M+1}}b_{M}(X_s^{M+1},\Lambda_s^{M+1},s)\d s\\
&\qquad+\int_0^{t\we\ta_{M}^{M+1}}\si_{M}(X_s^{M+1},\Lambda_s^{M+1},s)\d W_s,\\
\Lambda_{t\we\ta_M^{M+1}}^{M+1}&=i+\int_0^{t\we\ta_M^{M+1}}\int_0^\infty h(X_s^{M+1},\Lambda_{s-}^{M+1},z)N(\d z,\d s).
\end{align*}
By pathwise uniqueness of \eqref{equ_main_M},  $(X_{t\we\ta_M^M}^M,\Lambda_{t\we\ta_M^M}^M)=(X_{t\we\ta_M^{M+1}}^{M+1},\Lambda_{t\we\ta_M^{M+1}}^{M+1})$, for all $t\geq 0$, $\P$-a.s. Then we get that $\ta_{M+1}^{M+1}\geq\ta_M^{M+1}\geq\ta_M^M$. That means $\{\ta_k^k\}_{k\geq M,k\in\SS}$ is a sequence of stopping time and increasing. Let $\ta=\lim_{k\ra\infty}\ta_k^k$ and 
$$(X_t,\Lambda_t)=(X_t^k,\Lambda_t^k)~,t\leq \ta_k^k.$$
Then $\left\{\left(X_t,\Lambda_t\right)\right\}_{t\in[0,\ta)}$ is a local solution of \eqref{main_equ}. The uniqueness can be prove similarly as Theorem \ref{path_unique}.

\qed

We devote next theorem to the non-explosion of $\{(X_t,\Lambda_t)\}_{t\in[0,\ta)}$.  
\beg{thm}\label{non_explosion}
Fix $T>0$. $(X_t,\Lambda_t)$ is a local solution of \eqref{main_equ} starting from $(x,i)$. Assume that
\bequ\label{ineq_si}
\int_0^T\sup_{|y|+j\leq R}|\si(y,j,t)|_{HS}^2\d t<\infty,~R>0,\enqu
and there is a constant $\be>0$ such that the function
\bequ\label{equ_xj}
(y,j)\in\R^d\times \SS\mapsto \sum_{k=1}^\infty |k^\be-j^\be|q_{jk}(y)
\enqu
is locally bounded on $\R^d\times \SS$. If $b,~\si,~q_{ij}$ satisfy one of the following conditions, 
\beg{enumerate}
\item [$(1)$]  there exist a constant $p\geq 1$ and a non-negative function $C\in L^1([0,T],\R^+)$ such that
\beg{align*}
\ff {\sum_{k=1}^\infty\left(k^\be-j^\be\right)q_{jk}(y)} {\left(1+|y|^2\right)^{p}}&+\ff{2\<y,b(y,j,t)\>+(2p-1)|\si(y,j,t)|_{HS}^2} {1+|y|^2}\\
&\leq C(t)\left[1+\ff {j^\be} {\left(1+|y|^2\right)^p}\right],~(y,j)\in\R^d\times\SS,
\end{align*}
\item [$(2)$] there exist $\al\in(0,1]$, $c>0$ such that 
\beg{align}\label{ineq_exp}
&\ff{2\<y,b(y,j,t)\>+\left(1+2\al(1+|y|^2)^\al\right)|\si(y,j,t)|_{HS}^2} {1+|y|^2}\nonumber\\
&+\ff{\sum_{k\leq j}^\infty\left(k^\be-j^\be\right)q_{jk}(y)} {(1+|y|^2)^\al\exp\left\{(1+|y|^2)^{\al}\right\}}+\ff{\sum_{k> j}^\infty\left(k^\be-j^\be\right)q_{jk}(y)} {(1+|y|^2)^\al\exp\left\{e^{-\al cT}(1+|y|^2)^{\al}\right\}}\\
&\leq c\left[1+\ff {j^\be} {(1+|y|^2)^{\al}\exp\left\{e^{-\al cT}(1+|y|^2)^{\al}\right\}}\right],~(y,j)\in\R^d\times\SS,\nonumber
\end{align}
\end{enumerate}
then the solution $(X_t,\Lambda_t)$ is non-explosive on $[0,T]$.
\end{thm}
\beg{proof}
Let $\ta_M$ as in \eqref{ta_M}.
If $(1)$ holds, then by It\^o's formula, we have
\beg{align}\label{ineq_Xp}
&\d (1+|X_t|^2)^{p}\nonumber\\
&=p(1+|X_t|^2)^{p-1}\Big\{\Big[2\<X_t,b(X_t,\Lambda_t,t)\>+|\si(X_t,\Lambda_t,t)|_{HS}^2\nonumber\\
&\qquad+\ff {2(p-1)|\si^*(X_t,\Lambda_t,t)X_t|^2} {1+|X_t|^2}\Big]\d t+2\<X_t,\si(X_t,\Lambda_t,t)\d W_t\>\Big\}\\
&\leq p(1+|X_t|^2)^{p-1}\Big\{\Big[2\<X_t,b(X_t,\Lambda_t,t)\>+(2p-1)|\si(X_t,\Lambda_t,t)|_{HS}^2\Big]\d t\nonumber\\
&\qquad+2\<X_t,\si(X_t,\Lambda_t,t)\d W_t\>\Big\},~t<\ta_M.\nonumber
\end{align}
Moreover, by \eqref{ineq_si}
\beg{align*}
&\E\int_0^t(1+|X_{s\we\ta_M}|^2)^{p-1}|\si^*(X_{s\we\ta_M},\Lambda_{s\we\ta_M-},s)X_{s\we\ta_M}|^2\d s\\
&\leq M^2(1+M^2)^{p-1}\int_0^t\sup_{|x|+j\leq M}|\si(x,j,s)|_{HS}^2\d s\\
&<\infty.
\end{align*}
Hence, $\left\{\int_0^t(1+|X_{s}|^2)^{p-1}\<X_s,\si(X_{s},\Lambda_{s},s)\d W_s\>\right\}_{t\geq 0}$ is a local martingale. On the other hand, since 
$${\rm{supp}}\left(h(X_{t-},\Lambda_{t-},\cdot)\right)\subset \bigcup_{|x|\leq M, i\leq M,j\geq 1}\Ga_{i,j}(x)\subset \left[0,\sup_{|x|\leq M} \sum_{k=1}^{M+1}q_k(x)\right],~~t\leq\ta_M,$$
we have that $\{\Lambda_t\}_{t\in[0,\ta_M]}$ processes finite many jumps and
$$\Lambda_t=i+\int_0^t\int_{[0,K]}h(X_{s-},\Lambda_{s-},z)N(\d z,\d s),~t\leq \ta_M.$$
So, by It\^o's formula, we obtain that
\beg{align}\label{ineq_Lb}
\Lambda^\be_{t\we\ta_M}&=i^\be+\int_0^{t\we\ta_M}\int_{\R^+}\left[\left(\Lambda_{s-}+h(X_{s-},\Lambda_{s-},z)\right)^\be-\Lambda_{s-}^\be\right]N(\d z,\d s)\nonumber\\
&=i^\be+\int_0^{t\we\ta_M}\int_{\R^+}\left[\left(\Lambda_{s-}+h(X_{s-},\Lambda_{s-},z)\right)^\be-\Lambda_{s-}^\be\right]\tld N(\d z,\d s)\nonumber\\
&~\qquad+\int_0^{t\we\ta_M}\int_{\R^+}\left[\left(\Lambda_{s-}+h(X_{s-},\Lambda_{s-},z)\right)^\be-\Lambda_{s-}^\be\right]\d z\d s \nonumber\\
&=i^\be+\int_0^{t\we\ta_M}\int_{\R^+}\left[\left(\Lambda_{s-}+h(X_{s-},\Lambda_{s-},z)\right)^\be-\Lambda_{s-}^\be\right]\tld N(\d z,\d s)\\
&~\qquad+\int_0^{t\we\ta_M}\int_{\R^+}\sum_{k=1}^\infty\left[k^\be-\Lambda_{s-}^\be\right]\1_{[h(X_{s-},\Lambda_{s-},z)+\Lambda_{s-}=k]}(z)\d z\d s\nonumber\\
&=i^\be+\int_0^{t\we\ta_M}\int_{\R^+}\left[\left(\Lambda_{s-}+h(X_{s-},\Lambda_{s-},z)\right)^\be-\Lambda_{s-}^\be\right]\tld N(\d z,\d s)\nonumber\\
&~\qquad+\int_0^{t\we\ta_M}\sum_{k=1}^\infty\left[k^\be-\Lambda_{s-}^\be\right]q_{\Lambda_{s-}k}(X_{s-})\d s.\nonumber
\end{align}
Moreover, since \eqref{equ_xj} is locally bounded, we obtain that
\beg{align*}
\E \int_0^{t}\int_{\R^+}&\left|\left(\Lambda_{s\we\ta_M-}+h(X_{s\we\ta_M-},\Lambda_{s\we\ta_M-},z)\right)^\be-\Lambda_{s\we\ta_M-}^\be\right|\d z\d s\\
&\leq \E\int_0^t\sum_{k=1}^\infty\left|k^\be-\Lambda_{s-}^\be\right|q_{\Lambda_{s\we\ta_M-}k}(X_{s\we\ta_M-})\d s\\
&\leq \int_0^t\sup_{|x|+j\leq M}\sum_{k=1}^\infty\left|k^\be-j^\be\right|q_{jk}(x)\d s\\
&<\infty.
\end{align*}
Hence, $\left\{\int_0^{t}\int_{\R^+}\left[\left(\Lambda_{s\we\ta_M-}+h(X_{s\we\ta_M-},\Lambda_{s\we\ta_M-},z)\right)^\be-\Lambda_{s\we\ta_M-}^\be\right]\tld N(\d z,\d s)\right\}_{t\geq 0}$ is a  martingale and $\E\Lambda_{t\we\ta_M}^\be<\infty$. Then, following from \eqref{ineq_Xp} and \eqref{ineq_Lb}, we get that
\beg{align*}
\E \Big[(1+&|X_{t\we\ta_M}|^2)^{p}+p\Lambda_{t\we\ta_M}^\be \Big]-(1+|x|^2)^p-pi^\be\\
&\leq \E\int_0^{t\we\ta_M}\Big\{p(1+|X_s|^2)^{p-1}\Big[2\<X_s,b(X_s,\Lambda_s,s)\>+(2p-1)|\si(X_s,\Lambda_s,s)|_{HS}^2\Big]\\
&\qquad\qquad\qquad +\sum_{k=1}^\infty\left[k^\be-\Lambda_{s-}^\be\right]q_{\Lambda_{s-}k}(X_{s-})\Big\}\d s\\
&\leq  \E\int_0^{t\we\ta_M}C(s)\Big[(1+|X_{s-}|^2)^p+\Lambda_{s-}^\be\Big]\d s\\
&\leq \E\int_0^tC(s)\Big[(1+|X_{s\we\ta_M}|^2)^p+p\Lambda_{s\we\ta_M}^\be\Big]\d s,~t\leq T.
\end{align*}
By Gronwall's inequality, 
\beg{align*}
\E \Big[(1+&|X_{t\we\ta_M}|^2)^{p}+p\Lambda_{t\we\ta_M}^\be\Big]\leq \exp\left[\int_0^tC(s)\d s\right] \left\{(1+|x|^2)^p+pi^\be\right\},~t\leq T.
\end{align*}
Therefore, the solution to \eqref{main_equ} is non-explosive on $[0,T]$.

The proof of the claim under the condition $(2)$ is similar. Let $\la=\al c$. By It\^o's formula, we have that
\beg{align*}
&\d \exp\left\{e^{-\la t}(1+|X_t|^2)^\al\right\}\\
&\leq \al e^{-\la t}(1+|X_t|^2)^\al\exp\left\{e^{-\la t}(1+|X_t|^2)^\al\right\}\Big\{-\ff {\la} {\al}\\
&\quad+\ff{2\<x,b(x,j,t)\>+\left(1+2\al(1+|x|^2)^\al\right)|\si(x,j,t)|_{HS}^2} {1+|x|^2}\Big\}\d t\\
&\quad+2\al e^{-\la t}(1+|X_t|^2)^{\al-1}\<X_t,\si(X_t,\Lambda_t,t)\d W_t\>,~t<\ta_M,
\end{align*}
and
\beg{align*}
e^{-\la t\we\ta_M}\Lambda^\be_{t\we\ta_M}&=i^\be+\int_0^{t\we\ta_M}e^{-\la s}\int_{\R^+}\left[\left(\Lambda_{s-}+h(X_{s-},\Lambda_{s-},z)\right)^\be-\Lambda_{s-}^\be\right]\tld N(\d z,\d s)\nonumber\\
&~\quad+\int_0^{t\we\ta_M}e^{-\la s}\sum_{k=1}^\infty\left[k^\be-\Lambda_{s-}^\be\right]q_{\Lambda_{s-}k}(X_{s-})\d s-\la\int_0^{t\we\ta_M}e^{-\la s}\Lambda_s^\be\d s.
\end{align*}
Then, by \eqref{ineq_exp} and $\la =\al c$, we obtain that
\beg{align*}
&\E\left\{\exp\left\{e^{-\la t\we\ta_M}(1+|X_{t\we\ta_M}|^2)^\al\right\}+\al e^{-\la t\we\ta_M}\Lambda_{t\we\ta_M}^\be\right\}\\
&\leq \exp\{(1+|x|^2)^\al\}+\al i^\be-\al\E \int_0^{t\we\ta_M}e^{-\la s}\Lambda_{s-}^\be\d s\\
&\quad +\E\Big\{\int_0^{t\we\ta_M}\al e^{-\la s}(1+|X_s|^2)^\al\exp\left\{e^{-\la s}(1+|X_s|^2)^\al\right\}\Big[-\ff {\la} {\al}\\
&\quad+\ff{2\<X_s,b(X_s,\Lambda_{s-},s)\>+\left(1+2\al(1+|X_{s}|^2)^\al\right)|\si(X_{s},\Lambda_{s-},s)|_{HS}^2} {1+|X_{s}|^2}\\
&\quad+\ff{\sum_{k=1}^\infty\left(k^\be-\Lambda_{s-}^\be\right)q_{\Lambda_{s-}k}(X_s)} {(1+|X_s|^2)^\al\exp\left\{e^{-\la s}(1+|X_s|^2)^{\al}\right\}}\Big]\d s\Big\}\\
&\leq \exp\{(1+|x|^2)^\al\}+\al i^\be+(c-1)\al\E\int_0^{t\we\ta_M} e^{-\la s}\Lambda_{s-}^\be\d s.
\end{align*}
By Gronwall's inequality, there is $C>0$ which is independent of $M$ such that
$$\E\left\{\exp\left\{e^{-\la t\we\ta_M}(1+|X_{t\we\ta_M}|^2)^\al\right\}+\al e^{-\la t\we\ta_M}\Lambda_{t\we\ta_M}^\be\right\}\leq C.$$
Therefore, the proof is completed.

\end{proof}

\beg{rem}
The conditions in $(1)$ and $(2)$ are an extension of that used in \cite[Lemma 3.1]{XZ}, and they can be applied to the situation treated in \cite[Theorem 2.3]{Shao15}. In fact, in that case, we can let $\be=1$,~$p=1$. Comparing with \cite[Theorem 2.3]{Shao15}, our conditions allow that $\{q_i(\cdot)\}_{i\in\SS}$ are not necessary uniformly bounded on $\R^d$. However, it is a pity that  $q_i(x)$ is at most linear growth for the variable $i$.  
\end{rem} 
We present the following example to illustrate our conditions on $Q$-matrix $Q(x)$.
\beg{exa}
Let $p\geq 1$, $\ga>2$, $q_{jk}(x)=\ff {j+|x|^p} {|k-j|^\ga},~j\neq k,~k,j\in\SS$, and $C=\sum_{k=1}^\infty\ff 1 {k^\ga}$. Then 
\beg{align*}
&q_j(x)= \left(\sum_{k\geq j+1}^\infty\ff 1 {|k-j|^\ga}+\sum_{1\leq k\leq j}\ff 1 {|k-j|^\ga}\right)(j+|x|^p)\leq 2C(j+|x|^p),\\
&q_j(x)\geq (j+|x|^p)\sum_{k\geq j+1}^\infty\ff 1 {|k-j|^\ga}\geq C(j+|x|^p),~~j\in\SS.
\end{align*}
For all $\be\in[1,\ga-1)$, there exist positive constants $C_\be$, $C_{\ga,\be}$ such that
\beg{align*}
\sum_{k=1}^\infty \left(k^\be-j^\be\right) q_{jk}(x)&\leq C_\be\sum_{k\geq j+1}\ff {|k-j|^{\be}+|k-j|j^{\be-1}} {|k-j|^\ga} (j+|x|^p)\\
&\leq C_{\be,\ga}( j^\be+|x|^p) ,~j\in\SS.
\end{align*}

\end{exa}

\section{Strong Feller Property}
In this section, we shall investigate the strong Feller property for $(X_t,\Lambda_t)$ following the idea used in Theorem \ref{path_unique}. Let $f$ be a bounded measurable function on $\R^d\times\SS$, $P_tf(x,i)=\E f(X^{x,i}_t,\Lambda^{x,i}_t)$, and   $P_t^Kf(x,i)=\E f(X_t^{K},\Lambda_t^K)$ with $\left(X_0^K,\Lambda_0^K\right)=(x,i)$. Denote the transition semigroup associated with \eqref{equ_i} by $P_{t_0,t+t_0}^i$($P_t^i$ if $t_0=0$), i.e.
$$P_{t_0,t+t_0}^i f(x)=\E f(x_t^{x,i,t_0}).$$
\beg{thm}\label{thm_SF}
Assume that {\bf(H)} holds, moreover for each $i\in\SS$, $q_i$ is locally bounded and the semigroup $\{P_t^i\}_{t>0}$ generated by \eqref{equ_i}($t_0=0$) is strong Feller. Let $(x,i)\in\R^d\times \SS$, $M\in\SS,~M>|x|+i$, 
$$\ta^{x,i}_M=\inf\left\{t\geq 0~\Big|~|X^{x,i}_t|+\Lambda^{x,i}_t\geq M\right\}.$$
Assume that for each $x\in\R^d$, there exits $\de>0$ such that
\bequ\label{app_t}
\lim_{M\ra\infty}\sup_{|y-x|\leq \de}\P\left(\ta^{y,i}_M\leq t\right)=0.
\enqu
Then $P_t$ is strong Feller for $t>0$. 
\end{thm}

The discussion here begins with  the solution to \eqref{equ_K} $(X_t^K,\Lambda_t^K)$. But we will study \eqref{equ_K} ``path by path" of $\{N_t^K\}_{t\geq 0}$.  Fix $T>0$, $K>0$. Let $\om_2$ be a path of the  process $N^K_t$. Then $\om_2$ has finite jumps on $[0,T]$ at most. Denote  that $\De \om_2(t)=\om_2(t)-\om_2(t-)$. Let $(x^{(\om_2)}_t,\la^{(\om_2)}_t)$ be the solution of 
\bequ
\beg{cases}
\d x_t^{(\om_2)}=b(x_t^{(\om_2)},\la^{(\om_2)}_t,t)\d t+\si(x_t^{(\om_2)},\la_t^{(\om_2)},t)\d W_t,~~x_0^{(\om_2)}=x,\\
\la_t^{(\om_2)}=\sum_{s\leq t}h(x_s^{(\om_2)},\la_s^{(\om_2)},\De \om_2(s)),~~\la_0^{(\om_2)}=i.
\end{cases}
\enqu
Let $t_n$ be the $n$-th  discontinuous point of $\om_2$, and
$$m(\om_2)=\sup_{n}\{t_n\leq T\}.$$
\beg{lem}\label{lem_markov}
Assume that {\bf(H)} holds. Let $f$ be a bounded Borel measurable function on $\R^d\times\SS$. For each $\om_2$, there exists $g_{m(\om_2)}^{(\om_2)}:\R^d\times \SS\ra \R$ which is bounded and independent of $(x,i)$ such that
$$\E f(x_T^{(\om_2)},\la_T^{(\om_2)})=P_{t_1\we T}^i \Big[g_m^{(\om_2)}(\cdot,i)\Big](x).$$
\end{lem}
\beg{proof}
We shall prove this lemma by induction. For $0<T<t_1$, we have
$$\E f(x_T^{(\om_2)},\lambda_T^{(\om_2)})=\E f(x_T^{(\om_2)},i)=\P^i_T f(x).$$
If $T=t_1$, then
$$\E f(x_T^{(\om_2)},\la_T^{(\om_2)})=\E f(x_T^{(\om_2)},i+h(x_T^{(\om_2)},i,\De w_2(T)))=P_{T}^i\Big[ g^{(\om_2)}_1(\cdot,i)\Big](x),$$
where $g^{(\om_2)}(y,j)=f(y,j+h(y,j,\De w_2(t_1)))$, $(y,j)\in\R^d\times\SS$. So, on $[0,t_1]$, there exists $g^{(\om_2)}:\R^d\times \SS\ra\R$ which is bounded and independent of $(x,i)$ such that 
$$\E f(x_T^{(\om_2)},\lambda_T^{(\om_2)})=P^i_{t_1\we T}\left[g^{(\om_2)}(\cdot,i)\right](x),~T\in(0,t_1].$$
Next, we assume that for $T\in[0,t_n]$, there is $g^{(\om_2)}_n:\R^d\times \SS\ra \R$ which is independent of $(x,i)$ such that
$$\E f(x_T^{(\om_2)},\la_T^{(\om_2)})=P_{t_1\we T}^i \Big[g_n^{(\om_2)}(\cdot,i)\Big](x).$$
When $T\in (t_n,t_{n+1})$, due to the argument used in Lemma \ref{lem_path}, we have
\beg{align*}
\E f(x_T^{(\om_2)},\la_T^{(\om_2)})&=\E f(x_{T-t_n+t_n}^{(\om_2)},\la_{t_n}^{(\om_2)})=\E f\left(x^{(x^{(\om_2)}_{t_n},\la^{(\om_2)}_{t_n})}_{T-t_n},\la_{t_n}^{(\om_2)}\right)\\
&=\E \left\{\E \left[f\left(x^{(x^{(\om_2)}_{t_n},\la^{(\om_2)}_{t_n})}_{T-t_n},\la_{t_n}^{(\om_2)}\right)~\Big|~\left(x_{t_n}^{(\om_2)},\la_{t_n}^{(\om_2)}\right)\right]\right\}\\
&=\E \left\{\left[\left(P_{t_n,T}^{j}f(\cdot, j)\right)(y)\right]_{(y,j)=(x^{(\om_2)}_{t_n},\la^{(\om_2)}_{t_n})}\right\}\\
&\equiv \E\tld g_n^{(\om_2)}(x^{(\om_2)}_{t_n},\la^{(\om_2)}_{t_n}).
\end{align*}
If $T=t_{n+1}$, then 
\beg{align*}
\la^{(\om_2)}_{t_{n+1}}&=\la^{(\om_2)}_{t_n}+h(x^{(\om_2)}_{t_{n+1}},\la^{(\om_2)}_{t_n},\De \om_2(t_{n+1}))\\
&=\la^{(\om_2)}_{t_n}+h(x_{t_{n+1}-t_n}^{(x^{(\om_2)}_{t_n},\la^{(\om_2)}_{t_n})},\la^{(\om_2)}_{t_n},\De \om_2(t_{n+1})).
\end{align*}
Arguing as $T=t_1$, let 
$$\tld{\tld g}_n^{(\om_2)}(y,j)=P^j_{t_n,T}f(\cdot,j+h(\cdot,j,\De \om_2(t_{n+1})))(y).$$
Then
$$\E f(x_T^{(\om_2)},\la_T^{(\om_2)})=\E \tld{\tld g}_n^{(\om_2)}(x^{(\om_2)}_{t_n},\la^{(\om_2)}_{t_n}).$$
By induction, we prove the lemma.

\end{proof}

{\noindent\textbf{\emph{Proof of Theorem \ref{thm_SF}}}}\\
By \eqref{app_t}, $\left(X_t,\Lambda_t\right)$ is non-explosive. Let $K= \sup_{|x|\leq M} \sum_{k=1}^{M+1}q_k(x).$ According to \eqref{equ_MK}, $(X_t,\Lambda_t)=(X_t^K,\Lambda_t^K),~t<\ta^{x,i}_M$.
Let $\ta_n$ be the time of the $n$-th arrival of $N(K,\cdot)$, $m=\max\left\{n~|~\ta_n\leq t\right\}$. Then $\ta_n>0$, $\P$-a.s. Let $\E_{\P_2}$ be the distribution of the process $\{N^K_t\}_{t\in[0,T]}$ on the path space. We have
\beg{align*}
\E f(X^K_t,\Lambda^K_t)&=\E \left\{\E\left[ f(X_t^K,\Lambda^K_t)~\Big|N^K_\cdot\right]\right\}\\
&=\E_{\P_2} \left\{\E\left[ f(x_t^{(\om_2)},\la^{(\om_2)}_t)~\Big|N^K_\cdot=\om_2\right]\right\}\\
&=\E_{\P_2}P^i_{\ta_1\we t}g_{m(\om_2)}^{(\om_2)}(x).
\end{align*}
Then by strong Feller property of $P_t^i$ and the dominate convergence theorem, we obtain that
\beg{align}\label{equ_SF}
\lim_{y\ra x}P^K_t f(y,i)&=\lim_{y\ra x}\E f(X^K_t,\Lambda^K_t)=\lim_{y\ra x}\E_{\P_2}P^i_{\ta_1\we t}g_{m(\om_2)}^{(\om_2)}(x)\nonumber\\
&=\E_{\P_2}\lim_{y\ra x}P^i_{\ta_1\we t}g_{m(\om_2)}^{(\om_2)}(y)=\E_{\P_2}P^i_{\ta_1(\om_2)\we t}g_{m(\om_2)}^{(\om_2)}(x)\\
&=P_t^K f(x,i).\nonumber
\end{align}
On the other hand 
\beg{align*}
&\left| P_t f(x,i)-P_t^K f(x,i)\right|\leq \E\left| f(X_t,\Lambda_t)-f(X^K_t,\Lambda^K_t)\right|\\
&\leq \E\left\{\left| f(X_t,\Lambda_t)-f(X^K_t,\Lambda^K_t)\right|\1_{[\ta^{x,i}_M\leq t]}\right\}+\E\left\{\left| f(X_t,\Lambda_t)-f(X^K_t,\Lambda^K_t)\right|\1_{[\ta^{x,i}_M>t]}\right\}\\
&\leq 2|f|_{\infty}\P\left(\ta^{x,i}_M\leq t\right).
\end{align*}
So
\beg{align}\label{equ_sF}
&\left|P_t f(x,i)-P_t f(y,i)\right|\nonumber\\
&\leq \left| P_t f(x,i)-P_t^K f(x,i)\right|+\left| P_t f(x,i)-P_t^K f(x,i)\right|\nonumber\\
&\qquad+\left| P_t^K f(y,i)-P_t^K f(x,i)\right|\nonumber\\
&\leq 2|f|_{\infty}\left[\P\left(\ta^{x,i}_M\leq t\right)+\P\left(\ta^{y,i}_M\leq t\right)\right]+\left| P_t^K f(y,i)-P_t^K f(x,i)\right|.
\end{align}
Then it is stander to prove the strong Feller property of $P_t$ via \eqref{equ_sF}, \eqref{app_t}, \eqref{equ_SF} and Lemma \ref{lem_markov}.

\qed 

\beg{rem}
This theorem partly generalizes \cite[Theorem 3.2]{Shao15}, and works well when the diffusion in each environment is degenerated.
\end{rem}
Combining the results in the previous section, it is easy to get that
\beg{cor}
Fix $T>0$. Under the assumptions of Theorem \ref{path_unique} and Theorem \ref{non_explosion}, if for each $i\in\SS$, the semigroup $\{P_t^i\}_{t\in(0,T]}$ generated by \eqref{equ_i}($t_0=0$) is strong Feller holds, then $P_t$ has strong Feller property for $t>0$.
\end{cor}
\beg{proof}
We only have to prove \eqref{app_t}. Let $\ta_M^{x,i}$ as in Theorem \ref{thm_SF}. According to the proof of Theorem \ref{non_explosion}, there exists a locally bounded function $f:[0,T]\times \R^d\times \SS\ra [1,\infty)$ with  
$$\lim_{M\ra\infty}\left[\inf_{t\in[0,T],|y|+j\geq M}f(t,y,j)\right]=\infty$$
such that
$$\E f(t\we\ta_M,X_{t\we\ta_M},\Lambda_{t\we\ta_M})\leq C f(0,x,i),~(x,i)\in\R^d\times\SS$$
for some $C$ independent of $M,x,i$. Then
\beg{align*}
\P\left(t\geq \ta_M^{x,i}\right)&\leq \E \ff {f(t\we\ta_M,X_{t\we\ta_M},\Lambda_{t\we\ta_M})} {f(\ta_M,X_{\ta_M},\Lambda_{\ta_M})}\1_{[t\geq \ta_M^{x,i}]}\\
&\leq \ff {\E f(t\we\ta_M,X_{t\we\ta_M},\Lambda_{t\we\ta_M})\1_{[t\geq \ta_M^{x,i}]}} {\inf_{t\in[0,T], |y|+j\geq M}f(t,y,j)}\\
&\leq \ff {C f(0,x,i)} {\inf_{t\in[0,T], |y|+j\geq M}f(t,y,j)},~(x,i)\in\R^d\times\SS.
\end{align*}
Hence
$$\limsup_{M\ra\infty}\sup_{|y-x|\leq 1}\P\left(t\geq \ta_M^{y,i}\right)\leq \lim_{M\ra\infty}\sup_{|y-x|\leq 1}\ff {C f(0,y,i)} {\inf_{t\in[0,T], |y|+j\geq M}f(t,y,j)}=0.$$
\end{proof}
\bigskip

\noindent\textbf{Acknowledgements}


The author would like to thank Professor Feng-Yu Wang and Jinghai Shao for their useful suggestions.

\end{document}